\documentclass{article}
\usepackage{amsfonts}
\usepackage{amsmath}
\usepackage{color}

\newtheorem{theorem}{Theorem}

\newtheorem{corollary}[theorem]{Corollary}

\newtheorem{lemma}[theorem]{Lemma}

\newtheorem{remark}[theorem]{Remark}

\begin{document}

\title{Subcritical branching processes in random environment with
immigration stopped at zero\thanks{%
\noindent Doudou Li and Mei Zhang were supported by the Natural Science
Foundation of China under the grant 11871103, V.Vatutin was partially
supported by the High-End Foreign Experts Recruitment Program (No.
GDW20171100029).}}
\author{ Doudou Li\thanks{%
School of Mathematical Sciences \& Laboratory of Mathematics and Complex
Systems, Beijing Normal University, Beijing 100875, P.R. China. Email:
lidoudou@mail.bnu.edu.cn},~\ Vladimir Vatutin\thanks{%
Steklov Mathematical Institute, 8 Gubkin St., Moscow, 119991, Russia, Email:
vatutin@mi-ras.ru} \ and Mei Zhang\thanks{%
School of Mathematical Sciences \& Laboratory of Mathematics and Complex
Systems, Beijing Normal University, Beijing 100875, P.R. China. Email:
meizhang@bnu.edu.cn} }
\maketitle

\begin{abstract}
\noindent We consider subcritical branching processes with
immigration which evolve under the influence of a random environment and
study the tail distribution of life periods of such processes defined as the length of the time interval between the moment when first invader (or invaders) came to an empty site until the moment when the site becomes empty again. We prove that the tail
distribution decays with exponential rate. The main tools are the change of measure and some conditional
limit theorems for random walks.
\end{abstract}


\section{~ Introduction and statement of main results}

 Galton-Watson branching processes with immigration are among the popular models of branching processes.
 Different versions of such processes have found applications in physics, demography, biology and other fields of science.

 One of the problems being interesting from theoretical and practical points for such processes is the distribution of the so-called life period of a branching process with immigration defined as the length of the time interval between the moment when first invader (or invaders) came to an empty site until the moment when the site becomes empty again (see, for instance, \cite{4}, \cite{11}, \cite{14}, \cite{16}). Information on the length of such periods may be used, for example, in epidemiology, seismology and ecology. In the contents of epidemics, such periods correspond to the duration of outbreaks of diseases that do not lead to full epidemics and to the period of occupancy of sites in metapopulations \cite{8}. They may be used to analyse the waiting time interval for the end of earthquake aftershocks \cite{9} or plasmid incompatibility \cite{12}, \cite{13}, and for considering some other models of similar nature.

In this note we consider Galton-Watson branching processes allowing immigration and evolving in a
random environment. Individuals in such processes reproduce independently of
each other according to offspring distributions which vary in a random
manner from one generation to the other. In addition, a number of immigrants
join each generation independently of the development of the population and
according to the laws varying at random from generation to generation. A
formal definition of the process looks as follows. Let $\Delta =\left(
\Delta _{1},\Delta _{2}\right) $ be the space of all pairs of probability
measures on $\mathbb{N}_{0}=\{0,1,2,\ldots \}.$ Supplying $\Delta $ with the
component-wise metric of total variation we obtain a Polish space. Let $%
\mathbf{Q}=\{F,G\}$ be a two-dimensional random vector with independent
components
\begin{equation*}
F:=( F( \{ j\} ) ,j=0,1,...), \quad G:=(
G( \{ j\}) ,j=0,1,...)
\end{equation*}%
taking values in $\Delta $, and let $\mathbf{Q}_{n}=\{F_{n},G_{n}\},n=1,2,%
\ldots ,$ be a sequence of independent copies of $\mathbf{Q}$. The infinite
sequence $\mathcal{E}=\left\{ \mathbf{Q}_{1},\mathbf{Q}_{2},...\right\} $ is
called a random environment.

A sequence of $\mathbb{N}_{0}$-valued random variables $\mathbf{Y}=\left\{
Y_{n},\ n\in \mathbb{N}_{0}\right\} $ specified on the respective
probability space $(\Omega ,\mathcal{F},\mathbf{P})$ is called a branching
process with immigration in the random environment (BPIRE), if $Y_{0}$ is
independent of $\mathcal{E}$ and, given $\mathcal{E}$ the process $\mathbf{Y}
$ is a Markov chain with
\begin{equation*}
\mathcal{L}\left( Y_{n}|Y_{n-1}=y_{n-1},\mathcal{E}=\{\mathbf{q}_{1},\mathbf{%
q}_{2},...\}\right) =\mathcal{L}(\xi _{n1}+\ldots +\xi _{ny_{n-1}}+\eta _{n})
\end{equation*}%
for every $n\in \mathbb{N}:=\mathbb{N}_{0}\backslash \left\{ 0\right\} $, $%
y_{n-1}\in \mathbb{N}_{0}$ and $\mathbf{q}_{1}=\left( f_{1},g_{1}\right) ,%
\mathbf{q}_{2}=\left( f_{2},g_{2}\right) ,...\in \mathbf{Q}$, where $\xi
_{n1},\xi _{n2},\ldots $ are i.i.d. random variables with distribution $%
f_{n} $ and independent of the random variable $\eta _{n}$ with distribution
$g_{n} $. In the language of branching processes $Y_{n-1}$ is the $(n-1)$th
generation size of the population, $f_{n}$ is the distribution of the number
of children of an individual at generation $n-1$ and $g_{n}$ is the law of
the amount of immigrants joining generation $n$.

Along with the process $\mathbf{Y}$ we consider a branching process $\mathbf{%
Z}=\left\{ Z_{n},\ n\in \mathbb{N}_{0}\right\} $ in the random environment
(BPRE) $\mathcal{E}_{1}=\{ F_{1},F_{2},...\} $ which, given $%
\mathcal{E}_{1}$ is a Markov chain with $Z_{0}=1$ and, for $n\in \mathbb{N}$
\begin{equation*}
\mathcal{L}\left( Z_{n}|Z_{n-1}=z_{n-1},\mathcal{E}_{1}=(f_{1},f_{2},...)%
\right) =\mathcal{L}(\xi _{n1}+\ldots +\xi _{nz_{n-1}}).
\end{equation*}

An important role in studying BPRE and BPIRE is played by the so-called
associated random walk $\mathcal{S}=\left\{ S_{0},S_{1},...\right\} $. This
random walk has initial state $S_{0}$ and increments $X_{n}=S_{n}-S_{n-1}$, $%
n\geq 1$, defined as
\begin{equation*}
X_{n}:=\log m\left( F_{n}\right) .
\end{equation*}%
Here the increments $X_n$ are i.i.d. copies of the logarithmic mean offspring
number $X:=\log $ $m(F)$ with%
\begin{equation*}
m(F):=\sum_{j=0}^{\infty }jF(\{j\}) .
\end{equation*}%
We suppose that $X$ is a.s. finite.

We call a BPIRE $\mathbf{Y}$ supercritical if $\mathbf{E}[X]\in(0,+\infty],$
subcritical if $\mathbf{E}[X]\in[-\infty,0),$ and critical if either $%
\mathbf{E}[X] =0$ or $\mathbf{E}[X]$ does not exist.

It will be convenient to assume that if $Y_{n-1}=y_{n-1}>0$ is the
population size of the ($n-1)$th generation of $\mathbf{Y}$ then first $\xi
_{n1}+\ldots +\xi _{ny_{n-1}}$ individuals of the $n$th generation are born
and than $\eta _{n}$ immigrants join the population.

This agreement allows us to consider a modified version $\mathbf{W}=\left\{
W_{n},\ n\in \mathbb{N}_{0}\right\} $ of the process $\mathbf{Y}$ specified
as follows. Assume, without loss of generality, that $Y_{0}>0.$ Let $%
W_{0}=Y_{0}$ and for $n\geq 1$,
\begin{equation*}
W_{n}:=\left\{
\begin{array}{cc}
0, & \text{ if }T_{n}:=\xi _{n1}+\ldots +\xi _{nW_{n-1}}=0, \\[0.1in]
T_{n}+\eta _{n}, & \text{if }T_{n}>0.%
\end{array}%
\right.
\end{equation*}%
We call $\mathbf{W}$ a branching process with immigration stopped at zero
and evolving in the random environment $\mathcal{E}$.

The aim of the present paper is to study, under the annealed approach, the
tail distribution of the random variable
\begin{equation*}
\zeta :=\min \left\{ n\geq 1:W_{n}=0\right\}
\end{equation*}%
for subcritical BPIRE. Observe that the tail distribution of $\zeta $ for critical BPIRE's was described in \cite{6}.

With each pair of measures $(F,G)$ we associate the respective probability
generating functions
\begin{equation*}
F(s):=\sum_{j=0}^{\infty }F\left( \left\{ j\right\} \right) s^{j},\qquad
G(s):=\sum_{j=0}^{\infty }G\left( \left\{ j\right\} \right) s^{j}.
\end{equation*}%
Given the environment $\mathcal{E}=\left\{ (F_{n},G_{n}),n\in \mathbb{N}%
\right\} $, we construct the i.i.d. sequence of pairs of generating
functions
\begin{equation*}
F_{n}(s):=\sum_{j=0}^{\infty }F_{n}\left( \left\{ j\right\} \right)
s^{j},\qquad G_{n}(s):=\sum_{j=0}^{\infty }G_{n}\left( \left\{ j\right\}
\right) s^{j},\quad s\in \lbrack 0,1],
\end{equation*}%
and use below the convolutions of the generating functions $F_{1},...,F_{n}$
specified for $0\leq i\leq n-1$ by the equalities%
\begin{eqnarray*}
F_{i,n}(s):= &&F_{i+1}(F_{i+2}(\ldots F_{n}(s)\ldots )),\quad \\
F_{n,i}(s):= &&F_{n}(F_{n-1}(\ldots F_{i+1}(s)\ldots ))\ \text{ and }%
F_{n,n}(s):=s.
\end{eqnarray*}

We assume for convenience that $W_{0}=Y_{0}>0$ has the (random) probability
generating function%
\begin{equation*}
N(0;s):=\frac{G_{0}(s)-G_{0}(0)}{1-G_{0}(0)}
\end{equation*}%
where $G_{0}(s)\overset{d}{=}G(s)$. Other cases of initial distributions may
be considered in a similar way.

Denote
\begin{eqnarray*}
H_{n} &:&=\mathbf{E}\left[ \Big(1-G_{0}(F_{0,n+1}(0))\Big)%
\prod_{i=1}^{n}G_{i}(F_{i,n+1}(0))\right] , \\
H_{n}^{\ast } &:&=\mathbf{E}\left[ \frac{1-G_{0}(F_{0,n+1}(0))}{1-G_{0}(0)}%
\prod_{i=1}^{n}G_{i}(F_{i,n+1}(0))\right] , \\
R_{n} &:&=\mathbf{P}\left( \zeta >n\right) ,
\end{eqnarray*}%
and let%
\begin{equation*}
\mathcal{H}(s):=\sum_{n=0}^{\infty }H_{n}s^{n}\text{ , }\mathcal{H}^{\ast
}(s):=\sum_{n=1}^{\infty }H_{n}^{\ast }s^{n}\text{ and }\mathcal{R}%
(s):=\sum_{n=1}^{\infty }R_{n}s^{n}.
\end{equation*}

It is known (see, Lemma 1 in \cite{6}) that $\mathcal{R}(s)$ can be
calculated by the formula%
\begin{equation*}
\mathcal{R}(s)=\frac{s\mathcal{H}^{\ast }(s)+sR_{1}}{1-s\mathcal{H}(s)}.
\end{equation*}

The following restrictions are imposed on the distributions of $F$ and $G$.

\textbf{Hypothesis A1}. The BPRE is subcritical, i.e.
\begin{equation*}
-\infty\leq\mathbf{E}[X] <0
\end{equation*}
and either $-\infty<\mathbf{E}\left[ Xe^{X}\right] <0$ (the strongly subcritical case),
 or $\mathbf{E}\left[ Xe^{X}\right] =0$ (the intermediate subcritical case), or there is a number $0<\beta <1$ such that
\begin{equation*}
\mathbf{E}[Xe^{\beta X}]=0
\end{equation*}
(the weakly subcritical case).

Note that the subcritical BPRE's mentioned in Hypothesis A1 do not exhaust all possible cases of subcritical BPRE's. For instance, they do not include the subcritical BPRE's where $\mathbf{E}\left[ Xe^{tX}\right]=\infty$ for all $t>0$ (see \cite{15}) or where $\mathbf{E}\left[ Xe^{ tX}\right]<0$ for all $0\leq t\leq\beta$ with $\beta=\sup\{t\geq 0: \mathbf{E}\left[ e^{tX}\right]<\infty\}\in (0,1)$ (see \cite{5}).

One of the main tools in analyzing properties of BPRE and BPIRE is a change
of measure. We follow this approach and introduce a new measure $\mathbb{P}$
by setting, for any $n\in \mathbb{N}$ and any measurable bounded function $%
\psi :\Delta ^{n}\times \mathbb{N}_{0}^{n+1}\rightarrow \mathbb{R}$
\begin{equation}
\mathbb{E}[\psi (\mathbf{Q}_{1},\cdots ,\mathbf{Q}_{n},W_{0},\cdots, W_{n})]:=\gamma ^{-n}\mathbf{E}[\psi (\mathbf{Q}_{1},\cdots
,\mathbf{Q}_{n},W_{0},\cdots ,W_{n})e^{\delta S_{n}}],  \label{Cange_delta}
\end{equation}%
with
\begin{equation*}
\gamma :=\mathbf{E}[e^{\delta X}].
\end{equation*}%
Here $\delta =1$ for strongly and intermediate subcritical BPIRE and $\delta
=\beta $ for weakly subcritical BPIRE.

Observe that $\mathbf{E}[Xe^{\delta X}]=0$ translates into
\begin{equation*}
\mathbb{E}[X]=0.
\end{equation*}

We assume that under the new measure the following set of conditions holds
true.

\textbf{Hypothesis A2}. The distribution of $X$ is nonlattice. If a BPIRE\
is either intermediate or weakly subcritical then, with respect to $\mathbb{P}$, the distribution of $X$
belongs
to the domain of attraction of a two-sided stable law with index $%
\alpha \in (1,2]$.

Since $\mathbb{E}[X]=0$ for the intermediate or weakly subcritical  BPRE's, Hypothesis A2 provides existence of an increasing
sequence of positive numbers
\begin{equation}
a_{n}=n^{1/\alpha }l_{1}(n)  \label{Def_a}
\end{equation}%
with slowly varying sequence $l_{1}(1),l_{1}(2),...$ such that the
distribution law of $S_{n}/a_{n}$ converges weakly, as $n\rightarrow \infty $
to the mentioned two-sided stable law.

Our next assumption concerns the standardized truncated second moment of $F$,
\begin{equation*}
\vartheta (a):=\sum\limits_{j=a}^{\infty }j^{2}F(\{j\})/m(F)^{2},\quad a\in
\mathbb{N}.
\end{equation*}

Define $\log ^{+}x:=\log (\max (1,x))$.

\textbf{Hypothesis A3}.

1) If the BPRE is intermediate subcritical, then
\begin{equation*}
\mathbb{E}\left[ (\log ^{+}\vartheta (a))^{\alpha +\epsilon }\right] <\infty
\end{equation*}
for some $\epsilon >0$ and some $a\in \mathbb{N}$.

2) If the BPRE is strongly subcritical, then
\begin{equation*}
\mathbb{E}[\log ^{+}\vartheta (a)]<\infty
\end{equation*}
for some $a\in \mathbb{N}$.

Now we impose restrictions on the
immigration component.

\textbf{Hypothesis A4. }
\begin{equation*}
\mathbb{E}\left[ \frac{G_{0}^{\prime }(1)}{1-G_{0}(0)}\right] <\infty .
\end{equation*}

With Hypotheses A1-A4 in hands we are ready to formulate the main result of
this note.

\begin{theorem}
\label{T subcriticaltail}Let Hypotheses A1-A4 be satisfied. Then, as $%
n\rightarrow \infty $

1) if the equation $r\mathcal{H}(r)=1$ has a root $1<r<\gamma ^{-1}$, then%
\begin{equation*}
\mathbf{P}(\zeta >n)\sim \frac{r\mathcal{H}^{\ast }(r)+rR_{1}}{\mathcal{H}%
(r)+r\mathcal{H}^{\prime }(r)}r^{-n-1};
\end{equation*}

2) if the BPIRE is weakly subcritical and $\gamma ^{-1}\mathcal{H}(\gamma
^{-1})<1,$ then
\begin{equation}
\mathbf{P}(\zeta >n)\sim C\frac{\gamma ^{n}}{b_{n}},\ C\in \left( 0,\infty
\right)\label{New}
\end{equation}%
with $b_{n}:=na_{n}$;

3) if the BPIRE is weakly subcritical and $\gamma ^{-1}\mathcal{H}(\gamma
^{-1})=1,$ then%
\begin{equation*}
\mathbf{P}(\zeta >n)=o(\gamma ^{n}).
\end{equation*}
\end{theorem}

\begin{remark}
We show below that, under our conditions, the equation $r\mathcal{H}(r)=1$
always has a root $r\in (1,\gamma ^{-1})$ for strongly and intermediate
subcritical BPIRE.
\end{remark}

We note that Zubkov \cite{16}  considered a similar problem for a Galton-Watson
branching process with immigration $\{Y_{c}(n),n\geq 0\}$
evolving in the constant environment specified by probability generating functions $F(s)$ and $G(s), 0\leq s\leq 1$. He
investigated for this case the distribution of the so-called life period $\zeta _{c}$ initiated at time $N$ and defined as
\begin{equation*}
Y_{c}(N-1)=0,\, \min_{N\leq k<N+\zeta _{c}}Y_{c}(k)>0,\, Y_{c}(N+\zeta _{c})=0.
\end{equation*}
It was shown that if $G(0)>0$ and $F^{\prime}(1)<1$, i.e. the process $\{Y_{c}(n),n\geq 0\}$ is subcritical then $\mathbf{P}(\zeta >n)=(K+o(1))r_c ^{n}$ as $n\to\infty$ where $r_c\in (0,1)$  is an explicitly known constant and the parameter $K$ is (depending on some additional technical conditions) either positive or equal to zero. Thus, the form (\ref{New}) for the tail distribution of the random variable $\zeta$ in subcritical BPIRE's is different from those known for the ordinary subcritical Galton-Watson processes with immigration.

The distribution of life periods for other models of branching processes with immigration
evolving in a constant environment was analysed, for instance, in \cite{4}, \cite{11}, \cite{13} and \cite{14}.

Theorem \ref{T subcriticaltail} complements the main result of \cite{6} where it was shown that the tail distribution of $\zeta$ for a class of the critical BPIRE's stopped at zero behaves like $n^{-\kappa}l(n).$ Here $\kappa\in (0,1)$ and $l(n)$ is a function slowly varying at infinity.

In the sequel if no otherwise is stated, we write $h_{n}\sim Ck_{n},C>0$
if $\lim_{n\rightarrow \infty }h_{n}/k_{n}=C$, $h_{n}=O(k_{n})$ if $%
\limsup_{n\rightarrow \infty }h_{n}/k_{n}<\infty,$ and $h_{n}=o(k_{n})$ if $%
\lim_{n\rightarrow \infty }h_{n}/k_{n}=0.$ We also denote by $%
C,C_{1},... $ positive constants which may vary from place to place.

\section{\protect\bigskip Auxiliary results}

Our goal is to investigate the asymptotic properties of $H_{n}$ and $%
H_{n}^{\ast }$ and, having the asymptotics in hands, to find an asymptotic
representation for $R_{n}$ as $n\rightarrow \infty $. Observing that
\begin{equation*}
H_{n}=\gamma ^{n+1}\mathbb{E}\left[ (1-G_{0}(F_{0,n+1}(0)))%
\prod_{i=1}^{n}G_{i}\left( F_{i,n+1}(0)\right) e^{-\delta S_{n+1}}\right]
:=\gamma ^{n+1}\overline{H}_{n},
\end{equation*}%
we reduce the first problem to considering the asymptotic behavior of $%
\overline{H}_{n}$. Similar reduction may be performed for $H_{n}^{\ast }$.

Set
\begin{equation*}
M_{n}:=\max \left( S_{1},...,S_{n}\right), \quad L_{n}=:\min
(S_{0},S_{1},...,S_{n})
\end{equation*}%
and denote
\begin{equation*}
\tau (n):=\min \left\{ i\geq 0:S_{i}=L_{n}\right\} .
\end{equation*}

\begin{lemma}
\label{L upperweakly} Let Hypotheses A1-A2 be satisfied. If the process is
weakly subcritical, then for each $\varepsilon >0$, there exists $%
p=p(\varepsilon )$ such that%
\begin{equation*}
\mathbb{E}[(1-F_{0,n}(0))e^{-\beta S_{n}};\tau (n)\in \lbrack p,n-p]]<\frac{%
\varepsilon }{b_{n}}
\end{equation*}%
for all sufficiently large $n$.
\end{lemma}

\textbf{Proof} Note that%
\begin{eqnarray*}
1-F_{0,n}(0) &=&\mathbf{P}\left( Z_{n}>0|\mathcal{E}\right) \leq \min_{0\leq
k\leq n}\mathbf{P}\left( Z_{k}>0|\mathcal{E}\right) \\
&\leq &\min_{0\leq k\leq n}\mathbf{E}\left[ Z_{k}|\mathcal{E}\right]
=e^{\min_{0\leq k\leq n}S_{k}}.
\end{eqnarray*}%
Therefore, for each $p\in \left[ 1,n/2\right] $
\begin{equation*}
\mathbb{E}\left[ (1-F_{0,n}(0))e^{-\beta S_{n}};\tau (n)\in \lbrack p,n-p]%
\right] \leq \mathbb{E}\left[ e^{-\beta S_{n}}\cdot e^{\min_{0\leq k\leq
n}S_{k}};\tau (n)\in \lbrack p,n-p]\right] .
\end{equation*}%
We fix $k\in \lbrack p,n-p]$, set $S_{j}^{\ast }:=S_{k+j}-S_{k}$, $%
j=0,1,...,n-k$ and denote $L_{n-k}^{\ast }:=\min_{0\leq j\leq
n-k}S_{j}^{\ast }$. The duality property of random walks gives
\begin{eqnarray}
\mathbb{E}\left[ e^{-\beta S_{n}}\cdot e^{S_{k}};\tau (n)=k\right] &=&%
\mathbb{E}\left[ e^{(1-\beta )S_{k}}\cdot e^{-\beta S_{n-k}^{\ast }};\tau
(k)=k;L_{n-k}^{\ast }\geq 0\right]  \notag \\
&=&\mathbb{E}\left[ e^{(1-\beta )S_{k}};M_{k}<0%
\right] \mathbb{E}\left[ e^{-\beta S_{n-k}};L_{n-k}\geq 0\right] .
\label{RW_decomposition2}
\end{eqnarray}%
According to Proposition 2.1 in \cite{2}, for each $\theta >0$ there
exist positive constants $K_{i}=K_{i}(\theta ),i=1,2$ such that, as $n\rightarrow \infty $
\begin{equation}
\mathbb{E}\left[ e^{\theta S_{n}};M_{n}<0\right] \sim \frac{K_{1}}{b_{n}}%
,\quad \mathbb{E}\left[ e^{-\theta S_{n}};L_{n}\geq 0\right] \sim \frac{K_{2}%
}{b_{n}}.  \label{Asym_MaxMin}
\end{equation}%
We know by (\ref{Def_a}) that $b_{n}$ is a regularly varying sequence.
Therefore, for any $\varepsilon >0$ there exists an integer number $%
p=p(\varepsilon )$ such that
\begin{eqnarray*}
\mathbb{E}\left[ (1-F_{0,n}(0))e^{-\beta S_{n}};\tau (n)\in \lbrack p,n-p]%
\right] &\leq &C\sum\limits_{k=p}^{n-p}\frac{1}{b_{k}}\frac{1}{b_{n-k}} \\
&\leq &\frac{C_{1}}{b_{n}}\sum\limits_{k=p}^{\infty }\frac{1}{b_{k}}\leq
\frac{\varepsilon }{b_{n}}
\end{eqnarray*}%
for all sufficiently large $n$.

The lemma is proved.

\begin{lemma}
\label{L upperinter} Let Hypotheses A1-A2 be satisfied. If the process is
intermediate subcritical, then for each $\varepsilon >0$, there exists $%
p=p(\varepsilon )$ such that%
\begin{equation*}
\mathbb{E}[(1-F_{0,n}(0))e^{-S_{n}};\tau (n)\in \lbrack 0,n-p]]<\frac{%
\varepsilon }{n^{1-\alpha ^{-1}}l_{2}(n)},
\end{equation*}%
for all sufficiently large $n,$ where $l_{2}(1),l_{2}(2),...$ is a sequence slowly varying at infinity.
\end{lemma}

\textbf{Proof} It follows from Lemma 2.2 in \cite{3} that, as $%
n\rightarrow \infty $
\begin{equation*}
\mathbb{P}[M_{n}<0]\sim \frac{1}{n^{1-\alpha ^{-1}}l_{2}(n)}.
\end{equation*}%
Setting $\beta =1$ in (\ref{RW_decomposition2}) and using the arguments of
the preceding lemma we see that for any $\varepsilon >0$,
\begin{eqnarray*}
&&\mathbb{E}\left[ (1-F_{0,n}(0))e^{-S_{n}};\tau (n)\in \lbrack 0,n-p]\right]
\leq C_{3}\sum\limits_{k=1}^{n-p}\frac{1}{k^{1-\alpha ^{-1}}l_{2}(k)}\frac{1%
}{b_{n-k}} \\
&&\qquad\qquad\qquad\leq \frac{C_4}{b_{n}}\sum\limits_{1\leq k\leq n/2}%
\frac{1}{k^{1-\alpha ^{-1}}l_{2}(k)}+\frac{C_{4}}{n^{1-\alpha ^{-1}}l_{2}(n)%
}\sum\limits_{j=p}^{\infty }\frac{1}{b_{j}} \\
&&\qquad\qquad\qquad=O\left( \frac{1}{nl_{1}(n)l_{2}(n)}\right) +\frac{1}{%
n^{1-\alpha ^{-1}}l_{2}(n)}O\left( \frac{1}{p^{1/\alpha }l_{1}(p)}\right),
\end{eqnarray*}%
completing the proof.

\

To go further we need to perform two more changes of measure using the
right-continuous functions $U:\mathbb{R}$ $\rightarrow \lbrack 0,\infty )$
and $V:\mathbb{R}$ $\rightarrow [0,\infty)$ specified by%
\begin{equation*}
U(x):=1+\sum_{n=1}^{\infty }\mathbb{P}\left( S_{n}\geq -x,M_{n}<0\right)
,\quad x\geq 0,
\end{equation*}%
\begin{equation*}
V(x):=1+\sum_{n=1}^{\infty }\mathbb{P}\left( S_{n}<-x,L_{n}\geq 0\right)
,\quad x\leq 0.
\end{equation*}

It is known (see, for instance, \cite{1}\ and \cite{2}) that for any
oscillating random walk
\begin{equation}
\mathbb{E}\left[ U(x+X);X+x\geq 0\right] =U(x),\quad x\geq 0,  \label{Mes1}
\end{equation}%
\begin{equation}
\mathbb{E}\left[ V(x+X);X+x<0\right] =V(x),\quad x\leq 0.  \label{Mes2}
\end{equation}

Let $\mathcal{E}=\left\{ \mathbf{Q}_{1},\mathbf{Q}_{2},...\right\} $ be a
random environment and let $\mathcal{F}_{n},n\geq 1,$ be the $\sigma $-field
of events generated by the random vectors $\mathbf{Q}_{1},\mathbf{Q}_{2},...,%
\mathbf{Q}_{n}$ and the sequence $W_{0},W_{1},...,W_{n}$. The $\sigma $%
-fields $\{\mathcal{F}_{n},n\geq 1\}$  form a filtration $\mathfrak{F}$ and the increments $X_{n},n\geq 1$
of the random walk $\mathcal{S}$ are measurable with respect to the $\sigma $%
-field $\mathcal{F}_{n}$. We now introduce for each $n\geq 1$ a probability measure $\mathbb{P}_{(n)}^{+}$ on the $\sigma $-field $\mathcal{F}%
_{n}$ by means of the density
\begin{equation*}
d\mathbb{P}_{(n)}^{+}:=U(S_{n})I\left\{ L_{n}\geq 0\right\} d\mathbb{P}.
\end{equation*}%
In view of the martingale property (\ref{Mes1}) of $U$ the sequence of measures $\left\{
\mathbb{P}_{(n)}^{+},n\geq 1\right\} $ is consistent on the filtration $\mathfrak{F}$.
This and Kolmogorov's extension theorem show that we may assume without loss of generality that there
exists a probability measure $\mathbb{P}^{+}$ on $\mathfrak{F}$ such that
\begin{equation}
\mathbb{P}^{+}|\mathcal{F}_{n}=\mathbb{P}_{(n)}^{+},\ n\geq 1.
\label{DefMeasures}
\end{equation}

In the sequel we allow for arbitrary initial value $S_{0}=x$. Then, we write
$\mathbb{P}_{x}$ and $\mathbb{E}_{x}$ for the corresponding probability
measures and expectations. Thus, $\mathbb{P}=\mathbb{P}_{0}.$ Using this
agreement we rewrite (\ref{DefMeasures}) as
\begin{equation*}
\mathbb{E}_{x}^{+}\left[ O_{n}\right] :=\frac{1}{U(x)}\mathbb{E}_{x}\left[
O_{n}U(S_{n});L_{n}\geq 0\right] ,\ x\geq 0,
\end{equation*}%
for every $\mathcal{F}_{n}$-measurable random variable $O_{n}$. \

Similarly, the martingale property (\ref{Mes2}) of $V$  gives rise to probability measures $\mathbb{P}^{-}_{x},x\leq
0 $, and
\begin{equation*}
\mathbb{E}_{x}^{-}\left[ O_{n}\right] :=\frac{1}{V(x)}\mathbb{E}_{x}\left[
O_{n}V(S_{n});M_{n}< 0\right] ,\ x\leq 0.
\end{equation*}

We now come back to branching processes. To have a unified approach in
studying the asymptotic behavior of $H_{n}$ and $H_{n}^{\ast }$ as $%
n\rightarrow \infty $ we consider the sequence%
\begin{equation*}
B_{n}(s):=\mathbb{E}\left[ (1-B(F_{0,n}(s)))\prod%
\limits_{i=1}^{n}G_{i}(F_{i,n}(s))e^{-\delta S_{n}}\right] ,n\geq 1,
\end{equation*}%
where $B(s)$ is a (random) probability generating function which is
independent of the sequence $\mathbf{Q}_{n},n\geq 1$, and satisfies the
restriction

\textbf{Hypothesis A4*.}%
\begin{equation*}
\mathbb{E}\left[ B^{\prime }(1)\right] <\infty .
\end{equation*}

Taking $B(s)=G_{0}(s)$  and  $s=F_{n+1}(0)$ leads to $\gamma ^{-n-1}H_{n},$
while $B(s)=(G_{0}(s)-G_{0}(0))/\left( 1-G_{0}(0)\right) $ with the same $s$
gives $\gamma ^{-n-1}H_{n}^{\ast }$.

Our plan is to find asymptotic representations of $B_{n}(s)$ for all types
of subcritical BPIRE. To this aim we will use a decomposition
\begin{equation*}
B_{n}(s)=\sum_{k=0}^{n}B_{k,n}(s)
\end{equation*}%
where
\begin{equation*}
B_{k,n}(s):=\mathbb{E}\left[ (1-B(F_{0,n}(s)))\prod%
\limits_{i=1}^{n}G_{i}(F_{i,n}(s))e^{-\delta S_{n}};\tau (n)=k\right] .
\end{equation*}

\subsection{Weakly subcritical case}

In this subsection we prove the following statement.

\begin{theorem}
\label{T Hnweakly} Let Hypotheses A1-A2 and A4* be satisfied. If \ the
process is weakly subcritical with parameter $\beta \in (0,1),$ then for
each $s\in (0,1)$
\begin{equation*}
B_{n}(s)\sim \frac{C_{\beta }(s)}{b_{n}},\quad C_{\beta }(s)>0,
\end{equation*}%
as $n\rightarrow \infty $.
\end{theorem}

The idea of proving Theorem \ref{T Hnweakly} looks as follows. We show that,
for a fixed $k$ and $n\rightarrow \infty $%
\begin{equation*}
B_{k,n}(s)\sim C_{k}(s)\mathbb{E}\left[ e^{-\beta S_{n}};L_{n}\geq 0\right]
,\ B_{n-k,n}(s)\sim \hat{C}_{k}(s)\mathbb{E}\left[ e^{\left( 1-\beta \right)
S_{n}};\tau (n)=n\right]
\end{equation*}%
for some positive constants $C_{k}(s)$ and $\hat{C}_{k}(s)$, while $%
\sum_{k=p}^{n-p}B_{k,n}(s)$ is negligible in comparison with $1/b_{n}$ if $p$
is sufficiently large.

The proof of the asymptotic representations above is based on several
important statements established in \cite{2}. To check the applicability
of the statements we need to prove several preparatory lemmas.

Let $Z(k,n)$ be the number of particles at moment $n$ in a branching process
initiated at time $k$ by a single particle and $Z_{i}(k,n)$, $i=1,2,...$ be
independent probabilistic copies of $Z(k,n)$.

Put
\begin{equation*}
Y(k,n):=Z_{1}(k-1,n)+...+Z_{\eta _{k-1}}(k-1,n),\quad \Xi
(n):=\sum\limits_{k=1}^{n}Y(k,n)+\eta _{n},
\end{equation*}%
where we assume (with a slight abuse of notation) that $B(s)$ is the
probability generating function of $\eta _{0}$.

From now on we let $m:=\left[ n/2\right] ,$ where $\left[ x\right] $ stands
for the integer part of $x$ and write
\begin{eqnarray*}
B_{n}(s) &=&\mathbb{E}\left[ (1-B(F_{0,n}(s)))\prod%
\limits_{i=1}^{m}G_{i}(F_{i,n}(s))\prod%
\limits_{i=m+1}^{n}G_{i}(F_{i,n}(s))e^{-\beta S_{n}}\right] \\
&=&\mathbb{E}\left[ (1-F_{m,n}^{Y(1,m)}(s))F_{m,n}^{\Xi
(m)-Y(1,m)}(s)\prod\limits_{i=m+1}^{n}G_{i}(F_{i,n}(s))e^{-\beta S_{n}}%
\right] .
\end{eqnarray*}

Introduce two-dimensional random variables
\begin{equation*}
\mathcal{U}_{n}=\left( \mathcal{U}_{n1},\mathcal{U}_{n2}\right) :=\left(
e^{-S_{m}}Y(1,m),e^{-S_{m}}(\Xi (m)-Y(1,m))\right) ,
\end{equation*}%
\begin{equation}
\widetilde{\mathcal{V}}_{n}(s)=\left( \widetilde{\mathcal{V}}_{n1}(s),%
\widetilde{\mathcal{V}}_{n2}(s)\right) :=\left( (F_{m,n}(s))^{\exp
\{-(S_{n}-S_{m})\}},\prod\limits_{i=m+1}^{n}G_{i}(F_{i,n}(s))\right) ,
\label{tildeVn}
\end{equation}

and
\begin{equation*}
\mathcal{V}_{n}(s)=\left( \mathcal{V}_{n1}(s),\mathcal{V}_{n2}(s)\right)
:=\left( (F_{n-m,0}(s))^{\exp
\{-S_{n-m}\}},\prod\limits_{j=1}^{n-m}G_{j}(F_{j-1,0}(s))\right) .
\end{equation*}

\begin{lemma}
\label{L Unconvergence} If a BPIRE is weakly subcritical and Hypotheses A2
and A4* are valid then, for each $x\geq 0$
\begin{equation*}
\mathcal{U}_{n}\rightarrow \mathcal{U}_{\infty }:=\left( \mathcal{U}_{\infty
1},\mathcal{U}_{\infty 2}\right) \quad \mathbb{P}_{x}^{+}-a.s.
\end{equation*}%
as $n\rightarrow \infty $, where $\left( \mathcal{U}_{\infty 1},\mathcal{U}%
_{\infty 2}\right) $ is a random vector whose components are positive with
positive probabilities.
\end{lemma}

\textbf{Proof} Since the measure $\mathbb{P}_{x}^{+}$ imposes restriction
on the offspring probability laws of particles but not on the reproduction
of particles themselves, one can check that the random sequences
\begin{equation*}
e^{-S_{m}}Y(1,m),\quad e^{-S_{m}}(\Xi (m)-Y(1,m)),\ m=1,2,...
\end{equation*}%
form, correspondingly, a non-negative martingale and a submartingale with
respect to the filtration $\mathfrak{F}$. Hence, there exists a random
variable $\mathcal{U}_{\infty 1}$ such that, as $m\rightarrow \infty $
\begin{equation*}
e^{-S_{m}}Y(1,m)\rightarrow \mathcal{U}_{\infty 1},\quad \mathbb{P}%
_{x}^{+}-a.s.
\end{equation*}%
Since
\begin{equation*}
\mathbb{E}_{x}^{+}\left[ e^{-S_{m}}Y(1,m)\right] =e^{-x}\mathbb{E}\left[
B^{\prime }(1)\right] \in \left( 0,\infty \right) ,
\end{equation*}%
the random variable $\mathcal{U}_{\infty 1}$ is positive with a positive
probability.

Next, we claim that
\begin{equation*}
\sup\limits_{m}\mathbb{E}_{x}^{+}\left[ e^{-S_{m}}\Xi (m)\right] <\infty .
\end{equation*}%
If we prove this statement, then we may conclude that, as $m\rightarrow
\infty $
\begin{equation*}
e^{-S_{m}}(\Xi (m)-Y(1,m))\rightarrow \mathcal{U}_{\infty 2}\quad \mathbb{P}%
_{x}^{+}-a.s.
\end{equation*}%
where the random variable $\mathcal{U}_{\infty 2}$ is positive with a
positive probability in view of $e^{-S_{m}}(\Xi (m)-Y(1,m))\geq
e^{-S_{m}}Y(2,m).$

To establish the desired estimate recall that according to our change of
measure,
\begin{eqnarray*}
\mathbf{I}_{m}&:=&\mathbb{E}_{x}^{+}\left[ e^{-S_{m}}\Xi (m)\right] =\frac{1%
}{U(x)}\mathbb{E}_{x}\left[ e^{-S_{m}}\Xi (m)U(S_{m})I\left\{ L_{m}\geq
0\right\} \right] \\
&=&\frac{1}{U(x)}\sum\limits_{k=1}^{m}\mathbb{E}%
_{x}[e^{-S_{m}}Y(k,m)U(S_{m})I\left\{ L_{m}\geq 0\right\} ]+\frac{1}{U(x)}%
\mathbb{E}_{x}\left[ \eta _{m}e^{-S_{m}}U(S_{m})I\left\{ L_{m}\geq 0\right\} %
\right] .
\end{eqnarray*}%
Conditioning first on the environment $\mathcal{E}$ and then on $\eta
_{k-1},S_{0},S_{1},...S_{m}$ and observing that, for any $k<m$%
\begin{equation*}
\mathbb{E}_{x}\left[ e^{-S_{k}}U(S_{m})I\left\{ L_{m}\geq 0\right\} \right] =%
\mathbb{E}_{x}\left[ e^{-S_{k}}U(S_{k})I\left\{ L_{k}\geq 0\right\} \right]
\end{equation*}%
in view of (\ref{Mes1}), we obtain
\begin{eqnarray*}
\mathbf{I}_{m} &=&\frac{1}{U(x)}\sum\limits_{k=1}^{m}\mathbb{E}_{x}\left[
\eta _{k-1}e^{-S_{k-1}}U(S_{m})I\left\{ L_{m}\geq 0\right\} \right] +\frac{1%
}{U(x)}\mathbb{E}_{x}\left[ \eta _{m}e^{-S_{m}}U(S_{m})I\left\{ L_{m}\geq
0\right\} \right] \\
&\leq &\frac{\mathbb{E}\left[ \eta _{1}+\eta _{0}\right] }{U(x)}%
\sum\limits_{k=0}^{m}\mathbb{E}_{x}\left[ e^{-S_{k}}U(S_{m})I\left\{
L_{m}\geq 0\right\} \right] \\
&= &\frac{\mathbb{E}\left[ G^{\prime }(1)+B^{\prime
}(1)\right] }{U(x)}\sum\limits_{k=0}^{m}\mathbb{E}_{x}\left[
e^{-S_{k}}U(S_{k})I\left\{ L_{k}\geq 0\right\} \right] .
\end{eqnarray*}%
Since $U(y)$ is a renewal function, there exists a constant $C$ such that $%
U(y)\leq C(1+y)$ for all $y\geq 0$. Combining this estimate with the
inequality
\begin{equation*}
(1+y)e^{-y}\leq 2e^{-y/2},\quad y\geq 0,
\end{equation*}%
we see that
\begin{equation*}
\mathbf{I}_{m}\leq \frac{C}{U(x)}\sum\limits_{k=0}^{m}\mathbb{E}_{x}\left[
e^{-S_{k}}(1+S_{k})I\left\{ L_{k}\geq 0\right\} \right] \leq \frac{2C}{U(x)}%
\sum\limits_{k=0}^{m}\mathbb{E}_{x}\left[ e^{-S_{k}/2}I\left\{ L_{k}\geq
0\right\} \right] .
\end{equation*}%
Recalling (\ref{Asym_MaxMin}), it follows that, for all $m\in \mathbb{N}$
\begin{equation*}
\mathbf{I}_{m}\leq \frac{2C}{U(x)}\sum\limits_{k=0}^{\infty }\mathbb{E}_{x}%
\left[ e^{-S_{k}/2}I\left\{ L_{k}\geq 0\right\} \right] <\infty ,
\end{equation*}%
as desired.

The lemma is proved.

Denote%
\begin{equation*}
\Xi _{z}(n)=\sum_{j=1}^{z}Z_{j}(0,n),
\end{equation*}%
and introduce\ the random vector
\begin{equation*}
\mathcal{U}_{n}(z)=\left( \mathcal{U}_{n1}(z),\mathcal{U}_{n2}\right)
:=\left( e^{-S_{m}}\Xi _{z}(m),e^{-S_{m}}(\Xi (m)-Y(1,m))\right) .
\end{equation*}

Setting $B(s)=s^{z}$ we obtain the following statement.

\begin{corollary}
\label{L Unconvergence_z}Under the conditions of Lemma \ref{L Unconvergence}%
, for each $z\in \mathbb{N}$ and $x\geq 0$
\begin{equation*}
\mathcal{U}_{n}(z)\rightarrow \mathcal{U}_{\infty }(z):=\left( \mathcal{U}%
_{\infty 1}(z),\mathcal{U}_{\infty 2}\right) \quad \mathbb{P}_{x}^{+}-a.s.
\end{equation*}%
as $n\rightarrow \infty ,$ where $\left( \mathcal{U}_{\infty 1}(z),\mathcal{U%
}_{\infty 2}\right) $ is a random vector whose components are positive with
positive probabilities.
\end{corollary}

Now we deal with measure $\mathbb{P}^{-}$.

\begin{lemma}
\label{L Vnconvergence} If a BPIRE is weakly subcritical and Hypotheses A2
and A4* are valid then, for each fixed $s\in (0,1)$ and $x\leq 0$
\begin{equation*}
\mathcal{V}_{n}(s)=\left( \mathcal{V}_{n1}(s),\mathcal{V}_{n2}(s)\right)
\rightarrow \mathcal{V}_{\infty }(s):=\left( \mathcal{V}_{\infty 1}(s),%
\mathcal{V}_{\infty 2}(s)\right) \quad \mathbb{P}_{x}^{-}-a.s.
\end{equation*}%
as $n\rightarrow \infty $, where $\mathcal{V}_{\infty 1}(s)$ and $\mathcal{V}%
_{\infty 2}(s)$ are proper positive random variables.
\end{lemma}

\textbf{Proof} The fact that $\mathcal{V}_{n1}(s)\rightarrow \mathcal{V}%
_{\infty 1}(s)$ $\quad \mathbb{P}_{x}^{-}-a.s.$ as $n\rightarrow \infty $ is
a particular case of Lemma 3.2 in \cite{2}. To prove convergence of $%
\mathcal{V}_{n2}(s)$, note that given Hypotheses A2 and A4*,
\begin{equation*}
\sum\limits_{j=1}^{n-m}(1-G_{j}(F_{j-1,0}(s)))\leq \sum\limits_{j=1}^{\infty
}G_{j}^{\prime }(1)(1-F_{j-1,0}(0))\leq \sum\limits_{j=1}^{\infty
}G_{j}^{\prime }(1)e^{S_{j-1}}<\infty \quad \mathbb{P}_{x}^{-}-a.s.
\end{equation*}%
for every $s\in \lbrack 0,1]$. Hence, for each $s>0$
\begin{equation*}
\prod\limits_{j=1}^{n-m}G_{j}(F_{j-1,0}(s))\rightarrow \mathcal{V}_{\infty
2}(s):=\prod\limits_{j=1}^{\infty }G_{j}(F_{j-1,0}(s))>0\quad \mathbb{P}%
_{x}^{-}-a.s.
\end{equation*}

The lemma is proved.

For $u_{i}\geq 0,0\leq v_{i}\leq 1,t\geq 0\ (i=1,2)$ introduce the function
\begin{equation*}
\varphi (\mathbf{u},\mathbf{v},t)=\varphi
((u_{1},u_{2}),(v_{1},v_{2}),t):=v_{1}^{(u_{1}+u_{2})e^{t}}v_{2}.
\end{equation*}%
One may check that $\varphi $ is bounded and continuous within the specified
range of variables. For\ $z\in \mathbb{N}$, let
\begin{equation*}
J_{\nu }(s;z):=\int_{t\in (-\infty ,0]}\int_{\mathbf{u}\in R^{2}}\int_{%
\mathbf{v}\in R^{2}}\varphi (\mathbf{u},\mathbf{v},-t)\mathbb{P}^{+}(%
\mathcal{U}_{\infty }(z)\in d\mathbf{u})\mathbb{P}_{t}^{-}(\mathcal{V}%
_{\infty }(s)\in d\mathbf{v})\nu _{\beta }(dt),
\end{equation*}%
where
\begin{equation*}
\nu _{\beta }(dt):=K_{1}e^{\beta t}V(t)I\{t<0\}dt
\end{equation*}%
with scaling constant
\begin{equation*}
K_{1}^{-1}:=\int e^{\beta t}V(t)I\{t<0\}dt.
\end{equation*}

\begin{lemma}
\label{L_Positive}If a BPIRE is weakly subcritical and Hypotheses A2 and A4*
are valid then, for each $z\in \mathbb{N}$%
\begin{equation*}
\lim\limits_{n\rightarrow \infty }\frac{\mathbb{E}\left[ F_{0,n}^{z}(s)\prod%
\limits_{i=1}^{n}G_{i}(F_{i,n}(s))e^{-\beta S_{n}};L_{n}\geq 0\right] }{%
\mathbb{E}\left[ e^{-\beta S_{n}};L_{n}\geq 0\right] }=J_{\nu }(s;z).
\end{equation*}
\end{lemma}

\textbf{Proof} We write%
\begin{eqnarray*}
&&\mathbb{E}\left[ F_{0,n}^{z}(s)\prod\limits_{i=1}^{n}G_{i}(F_{i,n}(s))e^{-%
\beta S_{n}};L_{n}\geq 0\right] \\
&&\qquad =\mathbb{E}\left[ \left( F_{m,n}(s)\right) ^{\Xi (m)-Y(1,m)+\Xi
_{z}(m)}\prod\limits_{i=m+1}^{n}G_{i}(F_{i,n}(s))e^{-\beta S_{n}};L_{n}\geq 0%
\right] .
\end{eqnarray*}

Observing that
\begin{equation*}
\mathbb{E}\left[ F_{0,n}^{z}(s)\prod\limits_{i=1}^{n}G_{i}(F_{i,n}(s))e^{-%
\beta S_{n}};L_{n}\geq 0\right] =\mathbb{E}\left[ \varphi (\mathcal{U}%
_{n}(z),\widetilde{\mathcal{V}}_{n}(s),S_{n})e^{-\beta S_{n}};L_{n}\geq 0%
\right] ,
\end{equation*}%
where $\widetilde{\mathcal{V}}_{n}(s)$ is the same as in \eqref{tildeVn} and
using Theorem 2.7 in \cite{2} we complete the proof of the lemma.

\bigskip For $a>0$ and $u_{i}\geq 0,0\leq v_{i}\leq 1,z\leq 0\ (i=1,2)$ let
\begin{equation*}
\phi _{a}(\mathbf{u},\mathbf{v},z)=\phi
_{a}((u_{1},u_{2}),(v_{1},v_{2}),z):=\left( 1-v_{1}^{u_{1}e^{z}}\right)
v_{1}^{u_{2}e^{z}}v_{2}e^{-z}I\{z\geq -a\}.
\end{equation*}%
Clearly, $\phi _{a}$ is bounded and continuous in the specified domain. By
means of $\phi _{a}$ we specify, for $s\in \lbrack 0,1]$ the function
\begin{equation*}
J_{\mu }(s;a):=\int_{t\in \lbrack 0,\infty )}\int_{\mathbf{u}\in R^{2}}\int_{%
\mathbf{v}\in R^{2}}\phi _{a}(\mathbf{u},\mathbf{v},-t)\mathbb{P}_{t}^{+}(%
\mathcal{U}_{\infty }\in d\mathbf{u})\mathbb{P}^{-}(\mathcal{V}_{\infty
}(s)\in d\mathbf{v})\mu _{1-\beta}(dt),
\end{equation*}%
where
\begin{equation*}
\mu _{1-\beta}(dt):=K_{2}e^{-(1-\beta )t}U(t)I\{t\geq 0\}dt
\end{equation*}%
with scaling constant
\begin{equation*}
K_{2}^{-1}:=\int e^{-(1-\beta )t}U(t)I\{t\geq 0\}dt.
\end{equation*}

\begin{lemma}
\label{L taonnlimit} If a BPIRE is weakly subcritical and Hypotheses A2 and
A4* are valid then, for each $s\in \lbrack 0,1]$
\begin{equation*}
\lim\limits_{n\rightarrow \infty }\frac{B_{n,n}(s)}{\mathbb{E}\left[
e^{(1-\beta )S_{n}};\tau (n)=n\right] }=J_{\mu }(s;\infty ).
\end{equation*}
\end{lemma}

\textbf{Proof} We write
\begin{eqnarray*}
B_{n,n}(s) &=&\mathbb{E}\left[ \frac{1-B(F_{0,n}(s))}{e^{S_{n}}}%
\prod\limits_{i=1}^{n}G_{i}(F_{i,n}(s))e^{(1-\beta )S_{n}};\tau (n)=n\right]
\\
&=&\mathbb{E}\left[ \frac{1-B(F_{0,n}(s))}{e^{S_{n}}}\prod%
\limits_{i=1}^{n}G_{i}(F_{i,n}(s))e^{(1-\beta )S_{n}}I_{\{S_{n}<-a\}};\tau
(n)=n\right] \\
&&\qquad +\,\mathbb{E}\left[ \phi _{a}(\mathcal{U}_{n},\widetilde{\mathcal{V}%
}_{n}(s),S_{n})e^{(1-\beta )S_{n}};\tau (n)=n\right] \\
&&\qquad \qquad :=g(n;s;1)+g(n;s;2).
\end{eqnarray*}%
According (\ref{Asym_MaxMin})
\begin{eqnarray*}
g(n;s;1) &\leq &\mathbb{E}\left[ B^{\prime }(1)\frac{1-F_{0,n}(s)}{e^{S_{n}}}%
e^{(1-\beta )S_{n}}I\left\{ S_{n}<-a\right\} ;\tau (n)=n\right] \\
&\leq &\mathbb{E}[B^{\prime }(1)]\mathbb{E}\left[ e^{(1-\beta
)S_{n}}I\left\{ S_{n}<-a\right\} ;\tau (n)=n\right] \\
&\leq &\mathbb{E}[B^{\prime }(1)]e^{-(1-\beta )a/2}\mathbb{E}\left[
e^{(1-\beta )S_{n}/2};M_{n}<0\right] \\
&\leq & C\mathbb{E}[B^{\prime
}(1)]e^{-(1-\beta )a/2}/b_{n}
\end{eqnarray*}%
for all $n\in \mathbb{N}$ and $a>0$.

Further, we know from Lemmas \ref{L Unconvergence}--\ref{L Vnconvergence}
that the conclusion of Theorem 2.8 in \cite{2} holds for $\phi _{a}(%
\mathbf{u},\mathbf{v},z)$, i.e.%
\begin{equation*}
\lim_{n\rightarrow \infty }\frac{\mathbb{E}\left[ \phi _{a}(\mathcal{U}_{n},%
\widetilde{\mathcal{V}}_{n}(s),S_{n})e^{(1-\beta )S_{n}};\tau (n)=n\right] }{%
\mathbb{E}\left[ e^{(1-\beta )S_{n}};\tau (n)=n\right] }=J_{\mu }(s;a).
\end{equation*}%
Hence, letting $a$ to infinity we prove the lemma.

\textbf{Proof of Theorem \ref{T Hnweakly}} Since
\begin{equation*}
(1-B(F_{0,n}(s)))\prod\limits_{i=1}^{n}G_{i}(F_{i,n}(s))\leq B^{\prime
}(1)(1-F_{0,n}(s)),
\end{equation*}%
it follows from Lemma \ref{L upperweakly} that for any $\varepsilon >0$
\begin{equation*}
\sum_{k=p}^{n-p}B_{k,n}(s)\leq \mathbb{E}\left[ B^{\prime }(1)\right]
\mathbb{E}\left[ (1-F_{0,n}(s))e^{-\beta S_{n}};\tau (n)\in \lbrack p,n-p]%
\right] \leq \frac{\varepsilon \mathbb{E}[B^{\prime }(1)]}{b_{n}}
\end{equation*}%
for all sufficiently large $n$ and $p=p(\varepsilon )$.

Further, for fixed $k\leq p$, we take the expectation with respect to the $%
\sigma -$algebra\ $\mathcal{F}_{k}$ and obtain
\begin{eqnarray*}
&&\mathbb{E}\left[ (1-B(F_{0,n}(s)))\prod%
\limits_{i=1}^{n}G_{i}(F_{i,n}(s))e^{-\beta S_{n}};\tau (n)=k\right] \\
&&\qquad \qquad \qquad =\mathbb{E}\left[ e^{-\beta S_{k}}\Theta
_{1}(n-k;Y(1,k),\Xi (k)); \tau (k)=k \right] ,
\end{eqnarray*}%
where%
\begin{eqnarray*}
\Theta _{1}(n;z_{1},z_{2})&:=&\mathbb{E}\left[ \left(
1-F_{0,n}^{z_{1}}(s)\right)
F_{0,n}^{z_{2}-z_{1}}(s)\prod\limits_{i=1}^{n}G_{i}(F_{i,n}(s))e^{-\beta
S_{n}};L_{n}\geq 0\right] \\
&=&\mathbb{E}\left[ F_{0,n}^{z_{2}-z_{1}}(s)\prod%
\limits_{i=1}^{n}G_{i}(F_{i,n}(s))e^{-\beta S_{n}};L_{n}\geq 0\right] \\
&-&\mathbb{E}\left[ F_{0,n}^{z_{2}}(s)\prod%
\limits_{i=1}^{n}G_{i}(F_{i,n}(s))e^{-\beta S_{n}};L_{n}\geq 0\right] .
\end{eqnarray*}

Using Lemma \ref{L_Positive}, applying the dominated convergence theorem and
recalling (\ref{Asym_MaxMin}), we conclude that
\begin{eqnarray}
&&\lim\limits_{n\rightarrow \infty }\frac{B_{k,n}(s)}{\mathbb{E}\left[
e^{-\beta S_{n}};L_{n}\geq 0\right] }  \notag \\
&&\quad :=\mathbb{E}\left[ e^{-\beta S_{k}}
\left( J_{\nu }(s;\Xi (k)-Y(1,k))-J_{\nu }(s;\Xi (k))\right); \tau (k)=k \right] .
\label{Lim_Right}
\end{eqnarray}

Finally, we fix $j\geq 0$ and consider the expectation
\begin{equation*}
B_{n-j,n}(s)=\mathbb{E}\left[ (1-B(F_{0,n}(s)))\prod%
\limits_{i=1}^{n}G_{i}(F_{i,n}(s))e^{-\beta S_{n}};\tau (n)=n-j\right] .
\end{equation*}%
Denote $\overline{S}_{k}:=S_{n-j+k}-S_{n-j}\ (k=1,...,j)$ and let $\overline{%
\mathcal{F}}_{j}$ be the $\sigma $-algebra generated by
$\mathbf{Q}_{n-j+1},...\mathbf{Q}_{n}$. Taking the internal expectation with
respect to $\overline{\mathcal{F}}_{j}$ and supplying the respective
variables with bars $^{-}$ we see that
\begin{equation*}
B_{n-j,n}(s)=\mathbb{E}\left[ e^{-\beta \overline{S}_{j}}B_{n-j,n-j}(%
\overline{F}_{0,j}(s))\prod\limits_{i=1}^{j}\overline{G}_{i}(\overline{F}%
_{i,j}(s)); \overline{L}_{j}\geq 0\right] .
\end{equation*}%
Using Lemma \ref{L taonnlimit} and the dominated convergence theorem we
conclude that
\begin{eqnarray}
&&\lim_{n\rightarrow \infty }\mathbb{E}\left[ e^{-\beta \overline{S}_{j}}%
\frac{B_{n-j,n-j}(\overline{F}_{0,j}(s))}{\mathbb{E}\left[ e^{(1-\beta
)S_{n}};\tau (n)=n\right] }\prod\limits_{i=1}^{j}\overline{G}_{i}(\overline{F%
}_{i,j}(s))\overline{L}_{j}\geq 0\right]  \notag \\
&&\qquad \qquad =\mathbb{E}\left[ e^{-\beta \overline{S}_{j}}J_{\mu }\left(
\overline{F}_{0,j}(s);\infty \right) \prod\limits_{i=1}^{j}\overline{G}_{i}(%
\overline{F}_{i,j}(s))\overline{L}_{j}\geq 0\right] .  \label{Lim_Left}
\end{eqnarray}%
Combining (\ref{Lim_Right})-(\ref{Lim_Left}) with Proposition 2.1 in \cite%
{2}, we complete the proof.

\subsection{Intermediate and strongly subcritical cases}

In this subsection we find the asymptotics of $B_{n}(s)$ for intermediate
and strongly subcritical BPIRE.

\begin{theorem}
\label{T Hninter}Let Hypotheses A1-A3 and A4* be satisfied. If the process
is intermediate subcritical then, as $n\rightarrow \infty $
\begin{equation*}
B_{n}(s)\sim \frac{C}{n^{1-\alpha ^{-1}}l_{2}(n)},\quad C>0.
\end{equation*}
\end{theorem}

\textbf{Proof} Recalling that
\begin{equation*}
B_{k,n}(s)=\mathbb{E}\left[ (1-B(F_{0,n}(s)))\prod%
\limits_{i=1}^{n}G_{i}(F_{i,n}(s))e^{-S_{n}};\tau (n)=k\right],
\end{equation*}%
we have, for fixed $j\geq 0$
\begin{equation*}
B_{n-j,n}(s)=\mathbb{E}\left[ e^{-\overline{S}_{j}}B_{n-j,n-j}(\overline{F}%
_{0,j}(s))\prod\limits_{i=1}^{j}\overline{G}_{i}(\overline{F}_{i,j}(s)); \overline{L}_{j}\geq 0\right] .
\end{equation*}%
Using the duality property of random walks we see that
\begin{equation*}
B_{n,n}(s)=\mathbb{E}\left[ (1-B(F_{n,0}(s)))\prod\limits_{i=0}^{n-1}%
\widehat{G}_{i}(F_{i,0}(s))e^{-S_{n}};M_{n}<0\right] ,
\end{equation*}%
where $\widehat{G}_{i}$ are independent probabilistic copies of $G_{i}$.

Then, for any $\kappa >0$,
\begin{eqnarray*}
B_{n,n}(s) &=&\mathbb{E}\left[ (1-B(F_{n,0}(s)))\prod\limits_{i=0}^{n-1}%
\widehat{G}_{i}(F_{i,0}(s))e^{-S_{n}}I\{B^{\prime }(1)\leq \kappa \};M_{n}<0%
\right] \\
&&\qquad +\,\mathbb{E}\left[ (1-B(F_{n,0}(s)))\prod\limits_{i=0}^{n-1}%
\widehat{G}_{i}(F_{i,0}(s))e^{-S_{n}}I\{B^{\prime }(1)>\kappa \};M_{n}<0%
\right] \\
&&\qquad \qquad :=h(n;s;1)+h(n;s;2).
\end{eqnarray*}%
First observe that, as $n\rightarrow \infty $
\begin{equation*}
\frac{1-B(F_{n,0}(s))}{1-F_{n,0}(s)}\rightarrow B^{\prime }(1)\quad \mathbb{P%
}^{-}-a.s.
\end{equation*}%
and, by monotonicity and Hypothesis A3
\begin{equation*}
\frac{1-F_{n,0}(s)}{e^{S_{n}}}\geq \left( \frac{1}{1-s}+\sum\limits_{j=1}^{%
\infty }\vartheta _{j}(1)e^{S_{j}}\right) ^{-1}>0\quad \mathbb{P}^{-}-a.s.
\end{equation*}%
where $\vartheta _{j},j=1,2,...$ are i.i.d. copies of $\vartheta $. Hence,
there exists a positive random variable $\Theta (s)$ such that, as $n\to\infty$
\begin{equation*}
\frac{1-F_{n,0}(s)}{e^{S_{n}}}\rightarrow \Theta (s)\quad \mathbb{P}^{-}-a.s.
\end{equation*}%
Using the arguments similar to those applied to prove Lemma \ref{L
Vnconvergence}, we conclude that,  as $n\to\infty$
\begin{equation*}
\xi _{n}(s):=\prod\limits_{i=0}^{n-1}\widehat{G}_{i}(F_{i,0}(s))\rightarrow
\xi _{\infty }(s):=\prod_{i=0}^{\infty }\widehat{G}_{i}(F_{i,0}(s))>0\quad
\mathbb{P}^{-}-a.s.
\end{equation*}

Moreover,
\begin{equation*}
(1-B(F_{n,0}(s)))\prod_{i=0}^{n-1}\widehat{G}_{i}(F_{i,0}(s))e^{-S_{n}}\leq
\frac{1-B(F_{n,0}(s))}{1-F_{n,0}(s)}\frac{1-F_{n,0}(s)}{e^{S_{n}}}\leq
B^{\prime }(1).
\end{equation*}%
Hence it follows that
\begin{equation*}
\frac{h(n;s;2)}{\mathbb{P}\left[ M_{n}<0\right] }\leq \frac{\mathbb{E}\left[
B^{\prime }(1)I\{B^{\prime }(1)>\kappa \};M_{n}<0\right] }{\mathbb{P}\left[
M_{n}<0\right] }=\mathbb{E}\left[ B^{\prime }(1)I\{B^{\prime }(1)>\kappa \}%
\right]
\end{equation*}%
for all $n$ and, according to Lemma 2.5 in \cite{1}, as $n\rightarrow
\infty $
\begin{equation*}
\frac{h(n;s;1)}{\mathbb{P}\left[ M_{n}<0\right] }\rightarrow \mathbb{E}^{-}%
\left[ \Theta (s)B^{\prime }(1)I\{B^{\prime }(1)\leq \kappa \}\xi _{\infty
}(s)\right] .
\end{equation*}%
By these estimates and the arguments similar to those applied to check the
validity of Lemma \ref{L taonnlimit}, we conclude that, as $n\rightarrow
\infty $
\begin{equation*}
\frac{B_{n,n}(s)}{\mathbb{P}\left[ M_{n}<0\right] }\rightarrow \mathbb{E}^{-}%
\left[ B^{\prime }(1)\Theta (s)\xi _{\infty }(s)\right] .
\end{equation*}%
Combining this result with Lemma \ref{L upperinter} completes the proof.

\begin{theorem}
\label{T Hnstrongly}Let Hypotheses A1-A3 and A4* be satisfied. If the BPIRE\
is strongly subcritical then, for each $s\in (0,1)$%
\begin{equation*}
B_{n}(s)\sim C(s)>0
\end{equation*}%
as $n\rightarrow \infty $.
\end{theorem}

\textbf{Proof} The proof is based on the transformed measure $\mathbb{P}$.
In this case the inequality $\mathbf{E}[Xe^{X}]<0$ translates into
\begin{equation*}
\mathbb{E}[X]<0.
\end{equation*}%
Hence, the process is still subcritical under the probability measure $%
\mathbb{P}$ .

This fact, the equality
\begin{eqnarray*}
&&\mathbb{E}\left[ (1-B(F_{0,n}(s)))%
\prod_{i=1}^{n}G_{i}(F_{i,n}(s))e^{-S_{n}}\right] \\
&&\qquad \qquad =\mathbb{E}\left[ (1-B(F_{n,0}(s)))\prod_{i=0}^{n-1}\widehat{%
G}_{i}(F_{i,0}(s))e^{-S_{n}}\right] ,
\end{eqnarray*}%
the estimates
\begin{equation*}
(1-B(F_{n,0}(s)))\prod_{i=0}^{n-1}\widehat{G}_{i}(F_{i,0}(s))e^{-S_{n}}\leq
\frac{1-B(F_{n,0}(s))}{1-F_{n,0}(s)}\frac{1-F_{n,0}(s)}{e^{S_{n}}}\leq
B^{\prime }(1)
\end{equation*}%
and convergence
\begin{equation*}
(1-B(F_{n,0}(s)))\prod_{i=0}^{n-1}\widehat{G}_{i}(F_{i,0}(s))e^{-S_{n}}%
\rightarrow B^{\prime }(1)\Theta (s)\xi _{\infty }(s)\quad \mathbb{P}-a.s.
\end{equation*}%
as $n\rightarrow \infty, $ allow us to apply the dominated convergence
theorem to conclude that
\begin{equation*}
\lim\limits_{n\rightarrow \infty }B_{n}(s)=\mathbb{E}\left[ B^{\prime
}(1)\Theta (s)\xi _{\infty }(s)\right] >0.
\end{equation*}%
The theorem is proved.

\section{Proof of Theorem \protect\ref{T subcriticaltail}}

Our proof of Theorem \ref{T subcriticaltail} essentially uses the following
technical lemma.

\begin{lemma}
\label{L coefestimat}(see Theorem 1.4.6 in \cite{17}) Let
\begin{equation*}
\mathcal{T}(s)=\sum\limits_{n=0}^{\infty }\mathcal{T}_{n}s^{n}
\end{equation*}%
be a function with $\mathcal{T}_{n}\geq 0$ for all $n$. Assume that there
exist a number $\varrho >1$ and a function $l_{0}(n)$ slowly varying at
infinity such that
\begin{equation*}
\mathcal{T}_{n}\sim \frac{l_{0}(n)}{n^{\varrho }}
\end{equation*}%
as $n\rightarrow \infty $. If $\mathcal{C}(t)$ is an analytical function in
a domain containing the circle
\begin{equation*}
|t|\leq \mathcal{T}(1)=\sum\limits_{n=0}^{\infty }\mathcal{T}_{n},
\end{equation*}%
then
\begin{equation*}
\mathcal{C}(\mathcal{T}(s))=\sum_{n=0}^{\infty }c_{n}s^{n},\quad
\sum\limits_{n=0}^{\infty }|c_{n}|<\infty
\end{equation*}%
and
\begin{equation*}
c_{n}\sim \mathcal{C}^{\prime }(\mathcal{T}(1))\mathcal{T}_{n}.
\end{equation*}
\end{lemma}

Now everything is ready for proving Theorem \ref{T subcriticaltail}. We know
that%
\begin{equation*}
\mathcal{R}(s)=\frac{s\mathcal{H}^{\ast }(s)+sR_{1}}{1-s\mathcal{H}(s)}.
\end{equation*}%
Note that according to Theorems \ref{T Hnweakly}, \ref{T Hninter}, \ref{T
Hnstrongly} and the change of measure (\ref{Cange_delta}) there is a
positive constant $C\left( \delta \right) $ such that
\begin{equation}
\lim_{n\rightarrow \infty }\frac{H_{n}^{\ast }}{H_{n}}=C\left( \delta
\right) .  \label{Gen1}
\end{equation}%
Besides, for $\gamma =\mathbf{E}\left[ e^{\delta X}\right] $%
\begin{equation}
H_{n}\sim \left\{
\begin{array}{cc}
\frac{C_{\beta }}{b_{n}}\gamma ^{n+1} & \text{if the process is weakly
subcritical,} \\
&  \\
\frac{C\gamma ^{n+1}}{n^{1-\alpha ^{-1}}l_{2}(n)} & \text{if the process is
intermediate subcritical,} \\
&  \\
C\gamma ^{n+1} & \text{if the process is strongly subcritical.}%
\end{array}%
\right.   \label{Gen2}
\end{equation}

\textbf{Proof of Point 1) }Note that $\mathcal{H}(\gamma
^{-1})=\infty $ for the strongly and intermediate subcritical cases.\ Hence
a solution of the equation $r\mathcal{H}(r)=1$ within the interval $%
(1,\gamma ^{-1})$ always exists for these cases and $\mathcal{H}^{\ast
}(r)<\infty $. The same is true if $\gamma ^{-1}\mathcal{H}(\gamma ^{-1})>1$
for the weakly subcritical case. Taking these facts into account and
recalling point 3) of Theorem 1 in (\cite{7}, XIII.10) we conclude that
under the conditions of point 1) of Theorem \ref{T subcriticaltail}, as $%
n\rightarrow \infty $
\begin{equation*}
\mathbf{P}(\zeta >n)\sim \frac{r\mathcal{H}^{\ast }(r)+rR_{1}}{\mathcal{H}%
(r)+r\mathcal{H}^{\prime }(r)}r^{-n-1}.
\end{equation*}

\textbf{Proof of Point 2) } Setting
\begin{equation*}
\ \mathcal{T}(s):=\frac{s}{\gamma }\mathcal{H}\left( \frac{s}{\gamma }%
\right) :=\sum\limits_{n=1}^{\infty }\mathcal{T}_{n}s^{n},
\end{equation*}%
we see that, as $n\to\infty$
\begin{equation*}
\mathcal{T}_{n}\sim \frac{C_{\beta }}{b_{n}}.
\end{equation*}%
If $\mathcal{T}(1)=\gamma ^{-1}\mathcal{H}(\gamma ^{-1})<1$ then, taking
\begin{equation*}
\mathcal{C}(t):=\frac{1}{1-t}
\end{equation*}%
in Lemma \ref{L coefestimat} and writing
\begin{equation*}
\mathcal{C}(\mathcal{T}(s))=\frac{1}{1-\mathcal{T}(s)}=\sum\limits_{n=0}^{%
\infty }c_{n}s^{n},
\end{equation*}%
we conclude that, as $n\to\infty$%
\begin{equation*}
c_{n}\sim \mathcal{C}^{\prime }(\mathcal{T}(1))\mathcal{T}_{n}\sim \frac{%
C_{\beta }}{\left( 1-\gamma ^{-1}\mathcal{H}(\gamma ^{-1})\right) ^{2}}\frac{%
1}{b_{n}}.
\end{equation*}

Observing that
\begin{equation*}
\mathcal{R}\left( \frac{s}{\gamma }\right) =\frac{\frac{s}{\gamma }\mathcal{H%
}^{\ast }(\frac{s}{\gamma })+\frac{s}{\gamma }R_{1}}{1-\frac{s}{\gamma }%
\mathcal{H}\left( \frac{s}{\gamma }\right) }=\left( \frac{s}{\gamma }%
\mathcal{H}^{\ast }\left( \frac{s}{\gamma }\right) +\frac{s}{\gamma }%
R_{1}\right) \mathcal{C}(\mathcal{T}(s)),
\end{equation*}%
and using (\ref{Gen1})--(\ref{Gen2})\ we deduce, after evident estimates that%
\begin{equation*}
\frac{R_{n}}{\gamma ^{n}}=\sum_{k=1}^{n-1}\frac{H_{k}^{\ast }}{\gamma ^{k+1}}%
c_{n-k-1}+\frac{R_{1}}{\gamma }c_{n-1}\sim \frac{C}{b_{n}}
\end{equation*}%
as $n\rightarrow \infty $ as desired.

\textbf{Proof of Point 3) }Assume that $\mathcal{T}(1)=\gamma ^{-1}%
\mathcal{H}(\gamma ^{-1})=1$. Then%
\begin{equation*}
\overline{\mathcal{R}}(s):=\mathcal{R}\left( \frac{s}{\gamma }\right) =\frac{%
\frac{s}{\gamma }\mathcal{H}^{\ast }(\frac{s}{\gamma })+\frac{s}{\gamma }%
R_{1}}{1-\mathcal{T}(s)}:=\frac{\mathcal{G}(s)}{1-\mathcal{T}(s)}.
\end{equation*}%
By (\ref{Gen1}) and (\ref{Gen2})\

\begin{equation*}
\mathcal{G}(1)<\infty ,\quad \mathcal{T}(1)=1,\quad \mathcal{T}^{\prime
}(1)=\infty .
\end{equation*}%
Hence, applying  to the recurrent sequence $\left\{ \gamma ^{-n}R_{n},n\geq
1\right\} $ point 2) of Theorem 1 in (\cite{17}, XIII.10), we conclude that
\begin{equation*}
\lim_{n\rightarrow \infty }\frac{R_{n}}{\gamma ^{n}}=\frac{\mathcal{G}(1)}{%
\mathcal{T}^{\prime }(1)}=0.
\end{equation*}%
Theorem \ref{T subcriticaltail} is proved.

\end{document}